\documentclass[12p E.Lt,psfig]{article} % NEW LATEX

\usepackage{url}
\usepackage{epsfig}
\usepackage{multirow}

\def\convinlaw{\stackrel{{\cal L}}{\Longrightarrow }}             
\def\convinp{\stackrel{P}{\longrightarrow }}

\def\tends{\rightarrow}

\newcommand{\text}{\mbox}

\renewcommand{\marginpar}[1]{}

\begin{document}

\newtheorem{theorem}{Theorem}[section]
\newtheorem{lemma}{Lemma}[section]
\newtheorem{proposition}{Proposition}[section]
\newtheorem{corollary}{Corollary}[section]
\newtheorem{Remark}{Remark}[section]
\newtheorem{Example}{Example}[section]
\newtheorem{algorithm}{Algorithm}[section]

\newcommand{\e}{\varepsilon}
\newcommand{\ex}{\mbox{E}}
\newcommand{\var}{\mbox{Var}}
\newcommand{\xast}{X^{\ast}}
\newcommand{\btheta}{\mbox{\boldmath  $\theta $}}
\newcommand{\btau}{\mbox{\boldmath  $\tau $}}
\newcommand{\bdelta}{\mbox{\boldmath  $\delta $}}
\def\tends{\rightarrow}

\newpage

\title{Scalable subsampling: \\ computation, aggregation and inference}
\author{    Dimitris N. Politis
 \\ Department of
Mathematics \\
and Halicioglu Data Science Institute \\
University of California, San Diego\\ La Jolla, CA 92093-0112, USA
   \\  dpolitis@ucsd.edu } 
%\date{}
 
\maketitle

 \begin{abstract}  
Subsampling is a general statistical method developed in the 1990s  aimed at
estimating the sampling distribution of a statistic $\hat \theta _n$ in order to conduct 
    nonparametric inference such as the 
construction of confidence intervals and hypothesis tests.
Subsampling has seen a resurgence  in the {\it Big Data} era  
where the standard, full-resample size bootstrap can be infeasible to compute.
Nevertheless, even choosing a single random subsample of size $b$ can be
computationally challenging with both $b$ and the sample size $n$ being very large. 
In the paper at hand, we show how a set of appropriately chosen, non-random 
subsamples can be used to conduct effective---and computationally feasible---distribution estimation
via subsampling. Further,  we show how the same set of 
subsamples can be used to  yield  a procedure for subsampling aggregation---also known as
subagging---that is scalable with big data.
 Interestingly, the scalable subagging estimator  can be tuned to have
 the same (or better) rate of convergence as compared to $\hat \theta _n$.
The paper is  concluded by showing how to conduct inference,
e.g., confidence intervals, based on the scalable subagging estimator
instead of the original $\hat \theta _n$.

\noindent{\bf Keywords:} % \hspace{0.1cm}  
Bagging, Big Data, bootstrap, confidence intervals, distributed inference,  subagging. 
\end{abstract}

\newpage

\section{Introduction}
  
Assume data $X_1,\ldots, X_n$ that are independent, identically distributed (i.i.d.) taking values in an arbitrary 
space. Often, this space will be ${\bf R}^d$ but other choices exist, e.g., it can be a function space, 
a space of networks, etc. A statistic $\hat \theta_{n}=T_n(X_1,\ldots, X_n)$ is employed to estimate a
parameter $  \theta$ associated with the common distribution of the data. Assume that 
 $  \theta$ takes values on a normed linear space ${\bf \Theta}$ with  norm denoted by 
$||\cdot ||$.  

 Let $J_n(x)=P\{\tau_n g(\hat \theta_{n}-\theta)\leq x\}$ where the rate of convergence 
$\tau_n$ diverges to infinity as $n$ increases. 
In general,  the function $g(\cdot )$ may be taken to be the  norm  
$||\cdot ||$ but other choices are possible. For example, if ${\bf \Theta}={\bf R}^p,$ then 
$g(\cdot )$ can taken to be the sup-norm on ${\bf R}^p$; knowing (or estimating) the quantiles of $J_n$
would then lead to simultaneous confidence intervals
and/or simultaneous hypothesis tests for all  $p$ coordinates of $\theta$.
Note that if  ${\bf \Theta}={\bf R},$ then $g(\cdot )$
can be taken to be the identity function, leading to one-sided inference. 
We will also assume:
\\

\noindent 
{\bf Assumption A. } {\it There exists a non-degenerate probability distribution $J$, such that 
$J_n(x)\to J(x)$ as $n\to \infty$ for all $x$ points at which $J(x)$ is continuous.}
\\

Subsampling is a general statistical method developed in the 1990s  aimed at
estimating the sampling distribution $J_n$ in order to conduct 
    nonparametric inference such as the 
construction of confidence intervals and hypothesis tests; 
 see Politis, Romano and Wolf (1999) and the references therein.
To describe it, 
let the subsample size $b$ be a positive integer less than $n$, and 
consider all the subsets of size $b$ of the sample  $X_1,\ldots, X_n$.
There are $Q=$$n \choose b$ such subsets that can be ordered in 
an arbitrary fashion and denoted by ${\cal B}_j $ for $  j=1,\ldots, Q$.

Compute the subsample statistics $\hat \theta_{b,j}= T_b({\cal B}_j) $ for $  j=1,\ldots, Q$, and
the subsampling distribution 
$$  { } L_{ n,b} (x)  = \frac{1}{Q}\sum_{i=1}^{Q}
               1 \{ \tau_b g 
(\hat \theta_{  b,i}- \hat \theta_{  n} )    \leq x\}   .$$
Under Assumption A and the additional conditions
\begin{equation}
 b \to \infty \ \ \mbox{as} \ \ n\to \infty  \ \ \mbox{but with} \ \ b/n \to 0
\label{eq.1}
\end{equation}
and
\begin{equation}
 \tau _b/\tau _n \to 0
\label{eq.2} 
\end{equation} 
it was shown that  $   L_{ n,b} (x) \convinp J(x) $ for all $x$ points of continuity of $J$, 
where $\convinp$ denotes convergence in probability; see
Theorems 3.1  and 3.2 of Politis and Romano (1994). 
Note that $  { } L_{ n,b} (x) \convinp J(x) $ implies that 
$  { } L_{ n,b} (x) - J_n(x) \convinp 0$, i.e.,   $  { } L_{ n,b} (x) $ can be used as a proxy 
for the unknown $ J_n(x)$ so as to conduct statistical inference based on $\hat \theta _n$.
Also note that eq. (\ref{eq.2}) follows immediately if the rate of convergence satisfies
\begin{equation}
  \tau _n = n^\alpha {\cal L}(n)  \ \ \mbox{for some} \ \  \alpha >0,
\label{eq.3} 
\end{equation}
and some slowly varying function ${\cal L}(\cdot)$.

It was also recognized early on that if $n$ is large,  it is not realistic to compute 
$\hat \theta_{b,j}  $ for all $  j=1,\ldots, Q \ $ since $Q$ can be astronomically % exorbitantly
 large.  For that reason, 
  Corollary 2.2 of Politis and Romano (1994) showed that a stochastic approximation to 
$   L_{ n,b} (x) $  can be used instead. The stochastic approximation relies on 
$B$ randomly chosen subsamples from the set $\{{\cal B}_j $, $  j=1,\ldots, Q\}$. Here, $B$ 
is the practitioner's choice; it should be large,
tending to infinity as $n\to \infty$, but it can be chosen small enough for the
stochastic approximation  to be computable.

Note that choosing a random subsample from the set $\{{\cal B}_j $, $  j=1,\ldots, Q\}$
amounts to sampling $b$ values {\it without} replacement from the set  $\{X_1,\ldots, X_n\}$. This 
brings analogies to the so-called bootstrap with smaller resample size
where the sampling is {\it with} replacement---cf. Bickel and Freedman (1981),  
Bretagnolle  (1983), and  Swanepoel   (1986); see also 
Politis and Romano (1993), and Bickel,   G\"otze and  van Zwet (1997), 
and Corollary 2.3.1 of Politis,   Romano,   and Wolf  (1999)
for the comparison of sampling with vs.~without  replacement.

Subsampling has seen a resurgence  in the {\it Big Data} era of the 21st century
where Efron's (1979) full-resample size bootstrap can be infeasible to compute; see 
e.g., Jordan (2013), Kleiner et al. (2014), and Sengupta et al. (2016).
Nevertheless, even choosing a single random subsample of size $b$ can be
computationally challenging with both $b$ and $n$ very large; 
approximate solutions such as Poisson sampling may be needed as
described in Bertail,   Chautru,   and Cl\'emencon  (2017).

In the next section we show how a set of appropriately chosen, non-random 
subsamples can be used to conduct effective---and computationally feasible---distribution estimation
via subsampling. Further, in Section \ref{sec.subagg} we show how the same set of 
subsamples can be used to  yield  a procedure for subsampling aggregation---also known as
subagging---that is scalable with big data in an attempt to remedy computability issues
discussed in Section \ref{sec.comp}. 
 Interestingly, the scalable subagging estimator  can be tuned to have
 the same (or better) rate of convergence as compared to $\hat \theta _n$.
The paper is  concluded by providing details on how to conduct inference,
e.g. confidence intervals, based on the scalable subagging estimator
instead of the original $\hat \theta _n$.

\section{Scalable subsampling distribution estimation}
\label{sec.distr}

Recall the set of all size $b$ subsamples 
 $\{{\cal B}_j $,  $  j=1,\ldots, Q\}$, and re-arrange it so that the first subsamples 
are obtained as blocks of consecutive data points, i.e., 
%Compute the subsample statistics $\hat \theta_{j,b}= T_b({\cal B}_j) $ for $  j=1,\ldots, Q$, and
 %As  in Politis and Romano (1994, Section 3.2), let
${\cal B}_j=(X_{(j-1)h+1},X_{(j-1)h+2}, \ldots, X_{(j-1)h+b})$.
% $j$th partial-overlap, block--subsample of size $b$ that can be extracted from the data sequence $ \{X_1,\ldots,  X_{n}\} $.
Recall that the block size $b$ is an integer in $[1,n$], 
and so is $h$; %for some positive constant $C_1$, the parameter
in particular, $h$ % is   an integer in $[1,C_1   b$] 
  controls the amount of 
overlap (or separation) between ${\cal B}_j$ and ${\cal B}_{j+1}$.
If $h=1$, then the overlap is the maximum possible;  if $h\sim 0.2 \ b$, then there is an approximate 80\% 
overlap between ${\cal B}_j$ and ${\cal B}_{j+1}$;
if $h=b$, then there is $no$ overlap between ${\cal B}_j$ and ${\cal B}_{j+1}$
but these two blocks are adjacent; 
finally, if $h\sim 1.2 \ b$, then there is a block of about $0.2 \ b$ data points from the data sequence $ X_1,\ldots,  X_{n}  $ that separate the blocks  ${\cal B}_j$ and ${\cal B}_{j+1}$.
Note that, in general, $b$ and $h$ are functions of $n$, %and $h$ is a function of $b$,
but these dependences will not be explicitly denoted.

The collection of all available   block--subsample of size $b$, is then 
$\{{\cal B}_j, \ j=1,\ldots, q\}$  
where $q= \lfloor (n - b)/h\rfloor +1$ and $\lfloor \cdot \rfloor$ denotes the integer part function.
We claim that this non-random collection is sufficient for effective and computationally
feasible subsampling distribution estimation. To see why, note that subsampling using the 
aforementioned block--subsamples has been found to be consistent in the setting
where the data sequence $ X_1,\ldots,  X_{n}  $
is a finite stretch of a strong-mixing, stationary time series; see e.g. Politis and Romano (1994, Section 3.2).
Since the i.i.d. case can be considered as a special case of a stationary time series, 
the claim follows. 

To elaborate, we define the subsample statistics $\hat \theta_{b,j}= T_b({\cal B}_j) $ for $  j=1,\ldots, q$
as before, and construct the 
 subsampling distribution as 
\begin{equation}
 { } L_{ n,b,h} (x)  = \frac{1}{q}\sum_{i=1}^{q}
               1 \{ \tau_{  b} g(\hat \theta_{  b,i}- \hat \theta_{  n} )   \leq x\}   .
\label{eq.subtsh}
\end{equation}

\begin{proposition} 
Assume Assumption A, and conditions (\ref{eq.1}) and (\ref{eq.3}).
Also assume that either $h=1$, or that $h$ satisfies 
\begin{equation}
h\sim c_1 \ b \ \ \mbox{for some constant} \ \  c_1>0 .
\label{eq.5}
\end{equation}
Then, $   L_{ n,b,h} (x) \convinp J(x) $ for all $x$ points of continuity of $J$.
\end{proposition}

The Proposition follows from Corollary 3.2 of Politis  and Romano  (1994)  who worked under the
assumption that  $ 1\leq  h\leq b$; the case
where $h>b$ ---but still with $h=O(b)$---  can be proven in a similar way. The essence of the argument  is that
\begin{equation}
 EL_{ n,b,h} (x)   \approx J_b(x) \to J(x) \ \ \mbox{as} \ \  b\to \infty 
\label{eq.5.9}
\end{equation}
where  $x$ a point of continuity of $J$. In addition,  
\begin{equation}
 Var (L_{ n,b,h} (x)) = O(b/n).
\label{eq.6}
\end{equation}
Eq.~(\ref{eq.5.9}) and (\ref{eq.6}) together with 
Chebyshev's inequality   imply   $   L_{ n,b,h} (x) \convinp J(x) $.

Note that the bound (\ref{eq.6}) holds true regardless of the choice of $h$, i.e., whether 
$h=1$ or $h$ satisfies condition (\ref{eq.5}); it is just the 
proportionality constant in $O(b/n)$ that becomes smaller (but is bounded below) as $h$ decreases.  
Therefore, for reasons of parsimony and computational tractability, we will not propose using 
full-overlap block-subsamples, i.e., the case $h=1$, in what follows.
 Instead we will work under condition (\ref{eq.5}), 
in which case $q=O(n/b$). 
Hence, assuming $\hat \theta_{  n}$ can be computed in $O(n^\zeta$) operations
(for some constant  $\zeta >0$),  the construction of $L_{ n,b,h} $ and its quantiles
has computational complexity 
  $O(n^\zeta) +O(nb^{\zeta -1})=O(n^\zeta)$ which is the same 
%  computational complexity
  as computing the statistic  $\hat \theta_{  n}$ itself;
therefore, we may call the construction of $L_{ n,b,h} $
under  condition (\ref{eq.5}) as being {\it scalable}.

In terms of practicalities, the recommendation is: (a) choose an appropriate block size $b$,
  e.g.  use estimation or calibration as described in Chapter 9 of Politis et al. (1999),
or adapt the   method of G\"otze and Rackauskas (2001), 
and Bickel and Sakov (2008) to the scalable computation of 
$L_{ n,b,h} (x)$ that is based on block-subsamples; and (b) 
choose $h$ small, i.e., choose $c_1$ small in condition (\ref{eq.5}),  to make 
$q$ large (so that the the bound in (\ref{eq.6}) is small), but also lying 
within   feasibily/computability constraints. 

To give an example, suppose that $n=200,000$, and $b$ is chosen to be 500.
Recall that  $L_{ n,b,h} (x)$ is a sum of $q$ summands, where
 $q= \lfloor (n - b)/h\rfloor +1$.  If one chooses $c_1=2$, then $q\approx 200$ which is 
  feasible. If  a larger $q$ is deemed computationally feasible, one  could take 
$c_1=1$ or even $c_1=1/2$ leading to $q$ equalling 400 or 800 respectively.

Implementation problems might ensue  only if  $\zeta >1$
and $n$ is huge, but in this case even  
  the computation of the original statistic
$\hat \theta_{  n}$ may be problematic;
 we will deal with this issue in the next section. 

\section{Computability issues}  
\label{sec.comp}

As discussed in the last section, the statistic $\hat \theta_{  n}$ will generally be computable in $O(n^\zeta$) operations. 
If $\zeta$ is small, then no issues incur. Unfortunately, examples abound with $\zeta >1$,
making the   computability of  $\hat \theta_{  n}$  questionable in the Big Data era.

Some examples are as follows: 
\begin{enumerate}
\item  The $X_i$ are univariate, and $\hat \theta_{  n}$ is the sample mean (or median) of $X_1,\ldots, X_n$.
Then, $\zeta =1$.

\item  The $X_i$ take values in ${\bf R}^d$, and $\hat \theta_{  n}$ is the sample mean  of $X_1,\ldots, X_n$.
Then, $\hat \theta_{  n}$ is computable in $O(dn $) operations. If $d$ is a constant, then 
 $\zeta =1$ as above. However, it may be that $d$ grows with $n$; if $d$ grows
linearly with $n$,  then  $\zeta =2.$ 

\item Suppose that $X_i=(Y_i, W_i)$ where $Y_i$ is the univariate response associated with 
regressor $W_i$ that  takes values in ${\bf R}^d$; this is the standard regression situation. 
If $d$ is large, then  LASSO regression can be employed; see
Tibshirani (1996). 
 A popular method to compute the LASSO has computational complexity 
 $O(d^3 + d^2 n)$  
as long as $d<n $; see Efron et al. (2004).
 If $d$ grows
linearly with $n$, e.g., when $d\sim n/2$,  then  $\zeta =3.$ 
 
\end{enumerate}
 
There is a growing body of work dealing with the possibility that 
  the sample size $n$ is so large that it may not  be feasible to compute
 $\hat \theta_{  n}$.  One branch of this literature is devoted to 
`optimal subsampling' whose meaning is a little different than the subsampling-based 
inference discussed so far. In a nutshell, if $\hat \theta_{  n}$ is not computable, one
can use just one of the subsample statistics $\hat \theta_{ b,i}$ to estimate $  \theta $. 
The question `which one to use' is addressed by the literature on
`optimal subsampling'; see Yao and Wang  (2021) for a  review.

The problem with this approach is that the practitioner is effectivelly throwing out
the majority of the data. A Divide-and-Conquer alternative is proposed in the next section.

\section{Scalable subsampling aggregation}
\label{sec.subagg}
   
\subsection{Computation}  
\label{sec.SSA1}

Ensemble methods  are a family of techniques that combine
several estimators (or predictors)   with the objective
of producing a new estimator or predictor with better statistical properties.
One of the most popular and widely used ensemble methods is bootstrap aggregating
(Breiman, 1996), also known as {\it bagging}; it is a resampling technique whose
main purpose is to reduce the variability of a given  estimator/predictor
rendering it more {\it stable}.

Subsample aggregation, also known as {\it subagging}, was proposed by B\"uhlmann and Yu (2002). 
In the context of the present paper, the subagging estimator is defined as
$$\bar \theta_{b, SA} = \frac{1}{Q} \sum_{i=1}^Q \hat \theta_{  b,i}    $$
for an appropriate choice of $b$.  
Here, and for the remainder of the paper, we will assume that $\theta$ is univariate, i.e., 
 ${\bf \Theta}={\bf R}.$

Under some
regularity conditions, B\"uhlmann and Yu (2002) showed that
$$E\bar \theta_{b, SA} = E\hat \theta_{b,1} ,  \ \ \mbox{and} \ \ 
Var (\bar \theta_{b, SA})  \leq (b/n) Var (\hat \theta_{b,1}).$$
Hence, if the Bias of $\hat \theta_{b,1}$ is tolerable, 
subagging yields a welcome variance reduction.

With Big Data it is, of course,  infeasible to compute (and average) all 
 $Q=$$n \choose b$ subsample statistics.
As already mentioned,   even choosing a random subset of these $Q$ subsamples 
presents computational difficulties. Note, however, that although 
$\bar \theta_{b, SA}$ is an average of $Q$ values,   the variance is reduced by dividing by $n/b$ not $Q$;
the reason is that the $Q$ subsamples have typically high overlap, and hence their associated
 subsample statistics are highly dependent. 

We can achieve the same variance reduction effect by using just the first $q$ subsamples
in the ordering described in Section \ref{sec.distr}, i.e., 
${\cal B}_j=(X_{(j-1)h+1},X_{(j-1)h+2}, \ldots, X_{(j-1)h+b})$ for $j=1,\ldots,q$.
For better variance reduction ---as well as computational savings--- we will choose $h\geq b$ in what follows.
We therefore define the scalable   subagging estimator as
$$\bar \theta_{b, n,SS} = \frac{1}{q} \sum_{i=1}^q \hat \theta_{  b,i}   . $$

\begin{proposition} \label{prop.4.1}
Assume $h\geq b$. % as well as condition  (\ref{eq.5}). 
Then, 
$$E\bar \theta_{b, n,SS} = E\hat \theta_{b,1} ,  \ \ \mbox{and} \ \ 
Var (\bar \theta_{b, n,SS})  = q^{-1} Var (\hat \theta_{b,1})$$
where
$q= \lfloor (n - b)/h\rfloor +1$.
\end{proposition}
The proof is immediate noting that $ \hat \theta_{b,1}, \ldots, \hat \theta_{b,q}$
are i.i.d., their independence being ensured by the assumption $h\geq b$.

 \begin{Remark} \rm
Note that $q=O(n/b) $ since $h\geq b$; hence, 
 $ \bar \theta_{b, n,SS}$ can be computed with $O(q b^\zeta)=O(n b^{\zeta-1}) $  
operations, a significant savings over the $O(n^\zeta)$ operations needed for $\hat \theta_{n}$.
Recall the discussion of Section \ref{sec.comp} on (i) identifying the `optimal subsample' 
denoted by ${\cal B}_{j^*}$, and
then (ii) using as your final estimator the subsample statistic  $\hat \theta_{b,j^*}$.
The computational cost of part (ii) is of course $O(  b^\zeta)$; if you add to this
the nontrivial computational cost of part (i) ---that may require  numerical optimization---
the overall complexity may well exceed the 
   $O(n b^{\zeta-1}) $ needed to compute the scalable subagging estimator
$ \bar \theta_{b, n,SS}$.
\end{Remark}

%   \begin{Remark} \rm To achieve better variance reduction, it pays to choose $h$ small,
%e.g.,  when  $h=b$, then $q\approx n/b $. If $h$ instead satisfies condition  (\ref{eq.5})
%for some $c_1>1$, then $q$ will be of exact order $n/b $.
% \end{Remark}

We can further consider whether  $ \bar \theta_{b, n,SS}$ can have a  rate of  
convergence that is comparable to that of  $\hat \theta_{n}$; this is the subject of the
next subsection.

\subsection{Rate of convergence and choice of $b$ for scalable subagging}  
\label{sec.SSA2}

If the computational complexity of $ \bar \theta_{b, n,SS}$ is $O(n b^{\zeta-1}) $,
what is to stop us from taking $b$ very small, even $b=1$, to make it $O(n  ) $?
The answer is the bias of $\hat \theta_{b,1}$ that figures in
Proposition \ref{prop.4.1}. If the bias of $\hat \theta_{b,1}$ is zero (or negligible), then
indeed $b$ could be chosen very small. 

%{sec.comp}
Consider the following three Bias Conditions (BC): 
\begin{itemize}
\item [(I)] Estimator $\hat \theta_{n}$ is exactly unbiased, as is the case with linear statistics; see 
Example 1.4.1 of Politis et al. (1999). Then, $b$ can be taken to equal one but, of course, 
scalable subagging is not needed here as the 
 computational complexity of $\hat \theta_{n}$ is $O(n  ) $ already. 

\item [(II)] Estimator $\hat \theta_{n}$ is asymptotically unbiased, and its bias is asymptotically
negligible even after multiplication by $\tau _n$; in other words, the limit law $J$ of Assumption A
is centered at zero.

\item [(III)] Estimator $\hat \theta_{n}$ is asymptotically unbiased but its bias does not vanish 
 after multiplication by $\tau _n$; in other words, the limit law $J$ of Assumption A
is centered at a nonzero value. 
\end{itemize}
 
\noindent 
Although the subsampling distribution estimator $L_{n,b,h}$ can work under 
all above eventualities ---including BC 
case (III)---,  to investigate the rate of convergence of
scalable subagging we will work under the assumption of BC case (II). 
Note that an estimator falling under BC case (III) could be analytically {\it debiased}
---by subtracting from it a consistent estimate of its bias---,  allowing the 
debiased estimator to be handled under BC case (II).

Hence, we formulate the following assumption:
\\

\noindent 
{\bf Assumption B. } {\it Assume that $E\hat \theta_{n}^2 <\infty$ for 
all $n$, and that there exist constants $\gamma > \alpha >0$, $C\in {\bf R}-\{0\}$,
and $ \sigma^2 >0$ such that 
$$n^{ \gamma} (E\hat \theta_{n}-\theta) \to C ,  \ \ \mbox{and} \ \ 
Var (\tau_n \hat \theta_{n}) \to       \sigma^2    \ \ \mbox{as} \ \  n \to \infty $$
where $\tau _n = n^\alpha $.
}
\\
 
\noindent In the above, we have simplified eq.~(\ref{eq.3}) by omitting the slowly varying 
function ${\cal L}(n)  $; however, %this simplifies the subsequent discussion but
  ${\cal L}(n)  $
could  be incorporated in the technical arguments that follow if/when so desired. 
Note that the assumption  $\gamma > \alpha $, implies that the bias of 
$\hat \theta_{n}$ is negligible even after multiplication by $\tau _n$, as in BC case (II) above.
%However,, Assumption B  excludes estimators that are exactly unbiased.  
%$\tau _n$ is therefore the rate of convergence of $\hat \theta_{n}$.

Recall that in Section \ref{sec.distr} we operated under the assumption that $h$ is of the
same order of magnitude as $b$. So far, this was not required in the discussion of scalable subagging;
indeed, Proposition~\ref{prop.4.1} only needed $h\geq b$. 
For the sake of argument, in what follows we will entertain the possibility that 
$h$ grows at least as fast as $b$, and maybe faster. 
Therefore, we will  assume:
\begin{equation}
b\sim c_2 n^\beta \ \ \mbox{and} \ \ 
h\sim c_3 n^\delta  \ \ \mbox{as} \ \ 
n\to \infty \ \ \mbox{for   positive constants} \ \  c_2,c_3,
\ \ \mbox{and constants} \ \ 0< \beta \leq  \delta <1 .
\label{eq.7}
\end{equation}
  
The following
proposition shows that  $\bar \theta_{b, n,SS} $ can be tuned to have
 the same (or better) rate of convergence as compared to $\hat \theta _n$.
 We will use  the notation $a_n={\Theta} (d_n)$ to denote `exact order', i.e., that
there exist constants $\underline{c}, \bar c$ satisfying $\underline{c} \cdot  \bar c >0$, and  such that 
 $\underline{c} d_n \leq a_n \leq \bar c d_n.$

\begin{proposition} \label{prop.4.2}
Assume Assumption B, eq.~(\ref{eq.7}) and $h\geq b$.
Also assume $\alpha \leq 1/2$.
Choose a value of $\beta$
satisfying  $\beta\geq  \alpha /\gamma$, and 
%\begin{equation}
 % \frac{ \alpha }{\gamma} \leq  \beta< \frac{1}{2 (\gamma-\alpha)}.
%\label{eq.betabounds}
%\end{equation}
choose a value of $\delta$ satisfying 
\begin{equation}
1-2\beta(\gamma-\alpha) \leq \delta \leq 1+2\alpha (\beta -1).
\label{eq.deltabounds}
\end{equation}
%recalling the requirement $0<\delta $ as well. 
Then: 

(i)  
 $MSE(\bar \theta_{b, n,SS}) = O(\tau_n^{-2})$
where MSE is short for Mean Squared Error. 

(ii) If $\delta$ is chosen to equal
 the upper bound of eq.~(\ref{eq.7}),  i.e., if  $\delta= 1+2\alpha (\beta -1)$, then
$MSE(\bar \theta_{b, n,SS}) ={\Theta} (\tau_n^{-2}) $; in this case, 
$\bar \theta_{b, n,SS} $ has the same rate of convergence as $\hat \theta _n$.
However, if $\delta < 1+2\alpha (\beta -1)$, then
$MSE(\bar \theta_{b, n,SS}) = o(\tau_n^{-2}),$
and $\bar \theta_{b, n,SS} $ has {\em faster} rate of convergence than $\hat \theta _n$.

(iii) If $\delta$ can be chosen to equal
 the lower bound of eq.~(\ref{eq.7}),  i.e., if  $\delta=1-2\beta(\gamma-\alpha)$, then both 
$[Bias(\bar \theta_{b, n,SS})]^2$ and  $\ {Var (\bar \theta_{b, n,SS})  }$ have the
 same order of magnitude, namely they are both
(asymptotically) proportional to $n^{-2\beta \gamma}$ with generally different 
proportionallity constants; in this case, $\bar \theta_{b, n,SS}$ falls under BC case (III). 
However, if $\delta >1-2\beta(\gamma-\alpha)$, then
$[Bias(\bar \theta_{b, n,SS})]^2=o\left( {Var (\bar \theta_{b, n,SS})  }\right)$, 
i.e., $\bar \theta_{b, n,SS}$ falls under BC case (II).
 
\end{proposition}

\noindent {\bf Proof}. 
First note that $\bar \theta_{b, n,SS} $ is the average of $ \hat \theta_{b,1}, \ldots, \hat \theta_{b,q}$
that are  i.i.d.~due to the fact that $h\geq b$. Hence, $E\bar \theta_{b, n,SS}=E\hat \theta_{b,1}$, and 
$Bias(\bar \theta_{b, n,SS}) \sim C b^{-\gamma}) = C n^{-\beta \gamma}   $
  by Assumption B which further implies
\begin{equation}
Var (\bar \theta_{b, n,SS})  = q^{-1} Var(\hat \theta_{b,1})
\sim % {\Theta}(\frac{1}{ q b^{2\alpha}}) =
 \frac{\sigma^2}{   n^{1-\delta+ 2\alpha \beta}} ,
\label{eq.Var}
\end{equation}
since $q= \lfloor (n - b)/h\rfloor +1    ={\Theta}( n^{1-\delta})$.
Therefore,
$$MSE(\bar \theta_{b, n,SS}) \sim C^2 n^{-2\beta \gamma}+\frac{\sigma^2}{   n^{1-\delta+ 2\alpha \beta}}.$$

Choose %$\beta \geq \frac{1-\delta}{2(\gamma-\alpha)}$, i.e., 
$\delta \geq 1-2\beta(\gamma-\alpha)$,  to ensure
$ 2\beta \gamma \geq 1-\delta+ 2\alpha \beta$, yielding
$$MSE(\bar \theta_{b, n,SS}) =     {\Theta}(\frac{1}{   n^{1-\delta+ 2\alpha \beta}}) .$$
Finally, to have $ \frac{1}{   n^{1-\delta+ 2\alpha \beta}}=O(\frac{1}{   n^{  2\alpha  }})$, 
we need $\delta \leq 1+2\alpha (\beta -1)$; with equality, the two rates are the same,
i.e., if $\delta = 1+2\alpha (\beta -1)$, then
$ \frac{1}{   n^{1-\delta+ 2\alpha \beta}}  ={\Theta}(\frac{1}{   n^{  2\alpha  }})$.
%Check: can you get $1+2\alpha (\beta -1) \geq 1-2\beta(\gamma-\alpha)$ so that choice of $\delta $ is possible?
%YES need $\beta \geq \alpha/\gamma$ !!!  
$\diamond$  
 
\vskip .1in
 
\begin{Remark}
\label{re.Choice of b}
 \rm {\bf (Choice of $b$.)}
Elaborating on part (ii) of  Proposition \ref{prop.4.2},
if we wanted to maximize the rate of convergence of $\bar \theta_{b, n,SS} $ we would operate 
at the minimum value of $\delta$,  i.e., $\delta=1-2\beta(\gamma-\alpha)$, if the latter is a possible value 
for $\delta$. Recall that  eq.~(\ref{eq.7})  postulates   $\delta >0$; so 
 in order for $\delta$ to be chosen to equal $1-2\beta(\gamma-\alpha)$,
 the latter must be positive which can be ensured only when 
 \begin{equation}
\beta< \frac{1}{2 (\gamma-\alpha)} .
\label{eq.betabound}
\end{equation}

 %which would then yield
%$$MSE(\bar \theta_{b, n,SS}) = {\Theta} (n^{-2\beta \gamma}).$$
%\begin{corollary}
%Under the assumptions of    Proposition \ref{prop.4.2},
%the fastest rate of of convergence of $\bar \theta_{b, n,SS} $
%is $n^{\beta \gamma }, $ which is attainable when $\delta=1-2\beta(\gamma-\alpha)$.
%\end{corollary}
 
 %Under the premises of the above Corollary,
We could then choose $\beta$ to further
minimize $MSE(\bar \theta_{b, n,SS})$ which would entail choosing $
\beta$ to be as large as possible;   since $\beta \leq  \delta $
from  eq.~(\ref{eq.7}), it follows that we would need to take 
 $\beta =  \delta $. Furthermore, solving the equation $\beta =   1-2\beta(\gamma-\alpha)$
implies that the optimal value of $\beta$ is 
\begin{equation}
\beta=\frac{1}{1+2(\gamma -\alpha)}
\end{equation}
which is a value satisfying eq.~(\ref{eq.betabound}), thus allowing 
$\delta$ to operate 
at the minimum value of eq.~(\ref{eq.deltabounds}).
\end{Remark}
 
\begin{corollary}
\label{co:4.1}
Under the assumptions of    Proposition \ref{prop.4.2},
choose a value of $\beta$ satisfying
\begin{equation}
\frac{1}{1+2 (\gamma-\alpha)} \leq \beta< \frac{1}{2 (\gamma-\alpha)} ,
\label{eq.betabounds}
\end{equation}
and let $\delta=\beta$. Then, 
$$MSE(\bar \theta_{b, n,SS}) =     {\Theta}(\frac{1}{   n^{1-\beta+ 2\alpha \beta}}) .$$

Furthermore: 

\noindent
(i) If $\beta>\frac{1}{1+2 (\gamma-\alpha)}$, then 
$[Bias (\bar \theta_{b, n,SS})]^2=o\left( {Var (\bar \theta_{b, n,SS})  }\right)$, 
i.e., $\bar \theta_{b, n,SS}$ falls under BC case (II).

\noindent
(ii) If $\beta=\frac{1}{1+2 (\gamma-\alpha)}$, then 
$[Bias (\bar \theta_{b, n,SS})]^2={\Theta}\left( {Var (\bar \theta_{b, n,SS})  }\right)$, 
i.e., $\bar \theta_{b, n,SS}$ falls under BC case (III).
In addition, the choice $\beta=\frac{1}{1+2 (\gamma-\alpha)}$ minimizes 
the $MSE(\bar \theta_{b, n,SS})$,  yielding
\begin{equation} 
MSE(\bar \theta_{b, n,SS}) = {\Theta}(n^{-2  \gamma /[1+2(\gamma -\alpha)]})
\end{equation}
in which case the optimized rate of   convergence of $\bar \theta_{b, n,SS} $
is $n^{\gamma /[1+2(\gamma -\alpha)]}$.
\end{corollary}
\noindent Note that letting $\beta$ equal the lower bound 
$\frac{1}{1+2 (\gamma-\alpha)}$ kills two birds with one stone: (a) optimizes the 
rate of  convergence of $\bar \theta_{b, n,SS} $, and (b) minimizes the computational
complexity in computing $\bar \theta_{b, n,SS} $ making it
$O(n b^{\zeta -1})= O(  n^{1+\beta(\zeta -1)})
= O\left(n^{1+\frac{(\zeta -1)}{1+2 (\gamma-\alpha)}} \right).$

\vskip .1in

\begin{Remark}
\label{re.Choice of h}
 \rm {\bf (Choice of $h$.)} It was already establishes that $ \delta $ should be 
taken equal to its allowed minimum value to optimize $\bar \theta_{b, n,SS} $. Furthermore,
looking at eq.~(\ref{eq.Var}), it is apparent that to  optimize $Var (\bar \theta_{b, n,SS}) $
we would need to maximize $q$, i.e., minimize $h$. It is therefore recommended ---when  
computationally feasible--- to choose $h$ at its allowed minimum, i.e., choose $h=b$.
\end{Remark}

\vskip .1in
We now discuss some examples:
\begin{enumerate}
\item {\bf Linear statistics.} Consider a linear statistic $\hat \theta _n$ that can be represented
as $\hat \theta _n=n^{-1}\sum_{i=1}^n G(X_i)$ for some appropriate function $G$.
Note that $\hat \theta _n$ estimates $  \theta =EG(X_1)$, and is exactly unbiased for that. 
As mentioned under  the description of BC case (I), subagging is not needed here
because $\hat \theta _n$ can be easily computed. Furthermore, this case can not really fit under
the premises of Proposition \ref{prop.4.2} since Assumption B does not hold; recall that
Assumption B implies that $\hat \theta _n$ has nonzero bias.  However, 
we can intuite what would happen in this case by pretending that 
 Assumption B  holds (approximately) with a huge value of $\gamma$. 
Letting $\gamma \to \infty$, the above discussion implies that the optimal $\beta$ and 
$\delta$ should tend to zero. Hence, one should take $b=h=1$, in which
case $ \bar \theta_{b, n,SS}$ reduces to the original statistic  $\hat \theta _n$.

\item {\bf Approximately linear statistics.}
   A    statistic $\hat \theta _n$ is called approximately linear if it 
 can be represented
as $\hat \theta _n=n^{-1}\sum_{i=1}^n G(X_i)+o_P(n^{-1/2}) $.
Then, $\hat \theta _n$ is $\sqrt{n}$--consistent for $  \theta =EG(X_1)$, i.e., 
$\alpha=1/2$. In many interesting examples, e.g., sample quantiles, $M$-estimators, etc., 
it may be verified that the bias of  $\hat \theta _n$  is of order $1/n$, i.e., $\gamma =1$. In such a case, 
the two bounds in eq.~(\ref {eq.deltabounds}) collapse, yielding a single choice for $\delta$. Further 
optimizing MSE as in Remark \ref{re.Choice of b} implies
$\bar \theta_{b, n,SS}$ is $\sqrt{n}$--consistent as well, with the choice 
$\beta =  \delta=1/2$, i.e., block size increasing proportionally to 
 $\sqrt{n}$.

\item {\bf Nonparametric function estimators.}
Consider the case where $  \theta $ represents the value of function $f$ at 
a point of interest; the function $f$ can be a probability density, a spectral density, 
regression function, etc. that should be estimated in a nonparametric setting;
see   Rosenblatt (1991) for a unified description of 
nonparametric function estimation.
 Let $\hat \theta _n$  denote a 
 kernel-smoothed estimator of $  \theta $,
and suppose 
that a nonnegative kernel is used. In this case, the MSE--optimal 
bandwidth is ${\Theta} (n^{-1/5})$. However, this bandwidth choice brings 
 $\hat \theta _n$   under the realm of BC case (III), and 
the premises of Proposition \ref{prop.4.2} do not apply.
As an experiment, consider a degree of {\it undersmoothing}
in constructing $\hat \theta _n$. To fix ideas,  suppose 
  the bandwidth is chosen to be 
${\Theta} (n^{-1/4})$ instead, yielding
$$Bias(\hat \theta _n)=O(\frac{1}{n^{1/2}}) \ \ \mbox{and} \ \ 
Var(\hat \theta _n)={\Theta}(\frac{1}{n^{3/4}}) ;
$$
in this case, $\alpha=3/8$ and $\gamma=1/2$.
According to Remark \ref{re.Choice of b} and Corollary
\ref{co:4.1}, 
$  \beta $ and $  \delta $  should be optimally chosen to both equal $0.8$;
hence,  
  the  rate of   convergence of $\bar \theta_{b, n,SS} $ becomes $n^{2/5}$.
Interestingly, this rate is not only faster than the rate of $\hat \theta _n$ that used the
suboptimal  bandwidth ${\Theta} (n^{-1/4})$;  the  $n^{2/5}$ rate is actually the 
fastest rate achievable by {\it any} kernel-smoothed estimator 
that uses a nonnegative kernel with its associated MSE--optimal 
bandwidth.  Nevertheless, $\bar \theta_{b, n,SS} $ can be computed 
faster than  $\hat \theta _n$, and is thus preferable. 

\end{enumerate}

\vskip .1in
 
\subsection{Inference beyond point estimation}  
\label{sec.SSA3}

Having established that $\bar \theta_{b, n,SS} $ is a consistent estimator
whose rate of convergence towards $  \theta  $ is fast (and sometimes 
faster than that of  $\hat \theta _n$), the question now is  how to conduct inference,
e.g., confidence intervals, hypothesis tests, etc. based on $\bar \theta_{b, n,SS} $.

Recall that the assumption $h\geq b$ implies that $ \hat \theta_{b,1}, \ldots, \hat \theta_{b,q}$
are i.i.d.~but these values are not an increasing sequence; to see that, note that 
the value of $ \hat \theta_{b,1}$ changes with $b$ (which increases with $n$).
Rather, $ \hat \theta_{b,1}, \ldots, \hat \theta_{b,q}$
can be thought as the $n$th row of a triangular array with
  i.i.d.~entries, and common distribution  
given (approximately, and after centering and standardizing) by $J_b$.

Since $\bar \theta_{b, n,SS} $ is the sample mean of $ \hat \theta_{b,1}, \ldots, \hat \theta_{b,q}$, 
a Central Limit Theorem for triangular arrays ---such as Theorem B.O.1 of Politis et al.  (1999)---
implies the following useful result. 
 
\begin{corollary}
\label{co:4.2}
Assume the premises of  Corollary \ref{co:4.1}. 
In addition, assume    that there exist positive numbers $\epsilon$ 
and $\Delta$ such that
 $E|\hat \theta_n|^{2+\epsilon} \leq \Delta <\infty$ for  any $n$. 
Then, 
\begin{equation}
 n^{\frac{-1+\beta- 2\alpha \beta}{2}} \left( \bar \theta_{b, n,SS} -  \theta \right)
\convinlaw N( C_\beta, \sigma^2)
\label{eq.CLT}
\end{equation}
where $\convinlaw $ denotes convergence in law as $n\to \infty$. 
 Furthermore: 

\noindent
(i) If $\beta>\frac{1}{1+2 (\gamma-\alpha)}$, then 
$C_\beta =0$. 

\noindent
(ii) If $\beta=\frac{1}{1+2 (\gamma-\alpha)}$, then 
$C_\beta =C$. 

\noindent In the above, the constants $C$ and $\sigma^2$ are as defined in 
Assumption B.
\end{corollary}

Hence, under case (i) of the above, to construct an estimate of the distribution of $\bar \theta_{b, n,SS}$ 
all that is needed is a consistent estimator of $\sigma^2$; this is easy to obtain as 
the sample variance of the $n$th row of the triangular array, i.e., letting
\begin{equation}
\hat \sigma^2 =  \frac{b^{2\alpha}}{q}\sum _{i=1}^q \left( \hat \theta_{b,i} - \bar \theta_{b, n,SS} \right)^2 .
\label{eq.sample variance}
\end{equation}
For example, an approximate 95\% confidence interval  for $\theta$
under case (i) would be: $$\bar \theta_{b, n,SS} \pm 1.96\  \hat \sigma \cdot 
 n^{\frac{ 1-\beta + 2\alpha \beta}{2}} .$$
 
Actually, case (ii) of the Corollary \ref{co:4.2}  is more interesting since it ensures the
fastest rate of convergence of $\bar \theta_{b, n,SS}$. In this case, the nontrivial asymptotic
bias of the distribution of $\bar \theta_{b, n,SS}$ presents some difficulties at first inspection. 
Nevertheless,   subsampling comes again at the rescue since  eq.~(\ref{eq.CLT})
shows that $\bar \theta_{b, n,SS}$ has a well-defined asymptotic distribution; the fact that the latter
is not centered at zero is immaterial.
In other words, $\bar \theta_{b, n,SS}$ satisfies Assumption A with 
$\bar \theta_{b, n,SS}$ in place of $\hat \theta_{n}$. Hence, the scalable subsampling 
construction of Section \ref{sec.distr} can be applied to yield a consistent 
estimate of the distribution of $\bar \theta_{b, n,SS}$.

\begin{Remark}
\label{re.Connections with distributed inference}
 \rm {\bf (Connections with distributed inference.)}
The case $h=b$, i.e., splitting the sample into $q$ non-overlapping parts,  
is closely related to the classical  notion of $q$-fold cross-validation, as well 
as the more recent  notion of Divide-and-Conquer (DaC) methods; see 
 Jordan (2013).  
To elaborate,  the scalable subsampling estimator $\bar \theta_{b, n,SS} $
has been studied before in the following  DaC distributed inference   contexts (all with $h=b$):
     U-statistics by  Lin  and Xi  (2010);
  generalized linear models by Chen and Xie (2014); 
 M-estimators by Zhang, Duchi and Wainwright (2013); and 
a certain class of symmetric statistics  (that includes
 L-statistics  and   smoothed functions of the sample mean)  by Chen  and  Peng~(2021).
Note also that Bradic  (2016) employed subagging using non-overlapping blocks of data,   
  and applied it   to variable selection in large-scale regression.
The current section % \ref{sec.subagg}
is meant to serve many purposes. One is to show that these ideas  are universally applicable
under minimal assumptions, such as Assumption B; 
for example, the asymptotic normality results of Chen  and  Peng~(2021) actually  follow from 
  our Corollary \ref{co:4.2} simply by checking its premises. 
Furthermore, it is important to note that  the usefulness of scalable subagging and DaC distributed inference
extends well beyond the realm of asymptotically linear, $\sqrt{n}$--consistent statistics
that have been considered so far; see our Section \ref{sec.SSA2} including the 
example on nonparametric function estimation. 
Finally,  the interplay of the tuning parameters $h$ and $b$ opens up interesting
possibilities, e.g., the possibility that  $\bar \theta_{b, n,SS} $ has a faster
rate of convergence than  $\hat \theta_n $ itself.

\end{Remark}

\vskip .1in

\subsection{Possibly dependent data}
All subsampling constructions in this paper, including 
the scalable  subsampling distribution 
$ L_{ n,b,h} (x)  $, and 
scalable subagging estimator $\bar \theta_{b, n,SS} $, 
remain valid if there is (weak) dependence in the data, i.e., if
$X_1,\ldots,X_n$ are a stretch of a strictly stationary, strong mixing time series. 
The reason is that the choice of block-subsamples described
in Section \ref{sec.distr} and used throughout the paper is actually the choice
that is recommended in order to subsample time series; see e.g. Politis and Romano (1994). 
Hence, all results in this paper remain true as stated in the case where $X_1,\ldots,X_n$
are weakly dependent  with the  exception of the results in Section 
\ref{sec.SSA3} that may require a  little tweak. 

To see why, note that if 
$h=b$, then $ \hat \theta_{b,1}, \ldots, \hat \theta_{b,q}$ will be exactly
independent only if the data $X_1,\ldots,X_n$  are 
independent. If  $X_1,\ldots,X_n$  are stationary and strong mixing, 
we can ensure that $ \hat \theta_{b,1}, \ldots, \hat \theta_{b,q}$  
are approximately independent if we require that $h-b \to \infty$, 
e.g., letting $h=b+\lfloor \sqrt{b} \rfloor$; this would  ensure  that 
blocks ${\cal B}_j$ and ${\cal B}_{j+1}$ are separated by about $\sqrt{b} \ $ data points,
rendering them  approximately independent as $b$   increases with $n$. 
We   summarize this discussion in a corollary that covers the possibility of dependent data 
 with exponentially decreasing strong mixing coefficients.  

 \begin{corollary}
\label{co:4.3}
Assume the premises of  Corollary \ref{co:4.2}.  If  $h-b \to \infty$,  
then the results of Corollary~\ref{co:4.2} remain true even when the 
data $X_1,\ldots,X_n$ are a stretch of a strictly stationary, strong mixing time series,
with exponentially decreasing strong mixing coefficients. 
\end{corollary}
\noindent
Note that condition $h-b \to \infty$ ensures that the large sample variance of  $\bar \theta_{b, n,SS}$
is $\sigma^2$ as given in eq.~(\ref{eq.CLT}); 
 the asymptotic normality of $\bar \theta_{b, n,SS}$ (with different variance)  would  remain true for any $h\geq b$.  

An analog of Corollary \ref{co:4.3} can also be formulated when the strong mixing coefficients  only  decay polynomially fast. In that case, however,  the $\epsilon$  appearing in the moment assumption of Corollary~\ref{co:4.2}
can not be any positive number, i.e., it can not be taken arbitrarily close to zero. Rather, the minimum 
value of $\epsilon$ allowed would be dictated by the polynomial rate of decay of strong mixing; see
the assumptions of Theorem B.O.1 of Politis et al.  (1999).

\vskip .1in
\noindent
{\bf Acknowledgement.}  
Many thanks are due to Ery Arias-Castro,  Jelena Bradic
and Yiren Wang  for several helpful suggestions. 
Research  partially supported by NSF grant DMS 19-14556.

\end{document}